\documentclass[a4paper, 12 pt]{article}

\usepackage[T2A]{fontenc}
\usepackage[cp1251]{inputenc}
\usepackage[english]{babel}
\usepackage[tbtags]{amsmath}
\usepackage{amsfonts,amssymb}

\usepackage{mathrsfs}

\usepackage{geometry} 
\begin{document}
\title{$\alpha $, $\beta $-expansions of the  Riordan matrices of the associated subgroup.}
 \author{E. Burlachenko}
 \date{}

 \maketitle
\begin{abstract}
We consider the group of the matrices $\left( 1,g\left( x \right) \right)$ isomorphic to the group of formal power series $g\left( x \right)=x+{{g}_{2}}{{x}^{2}}+...$  under composition: $\left( 1,{{g}_{2}}\left( x \right) \right)\left( 1,{{g}_{1}}\left( x \right) \right)=\left( 1,{{g}_{1}}\left( {{g}_{2}}\left( x \right) \right) \right)$. Denote $P_{k}^{\alpha }=\left( 1,x{{\left( 1-k\alpha {{x}^{k}} \right)}^{{-1}/{k}\;}} \right)$. Matrix $\left( 1,g\left( x \right) \right)$ is decomposed into an infinite product of the matrices $P_{k}^{\alpha }$ with suitable exponents in two ways  –  to left-handed and right-handed products with respect to the matrix $P_{1}^{{{\alpha }_{1}}={{\beta }_{1}}}$: $\left( 1,g\left( x \right) \right)=...P_{k}^{{{\alpha }_{k}}}...P_{2}^{{{\alpha }_{2}}}P_{1}^{{{\alpha }_{1}}}=P_{1}^{{{\beta }_{1}}}P_{2}^{{{\beta }_{2}}}...P_{k}^{{{\beta }_{k}}}...$. 
We obtain two formulas expressing the coefficients of the series ${{\left( {g\left( x \right)}/{x}\; \right)}^{z}}$ in terms of the expansion coefficients ${{\alpha }_{i}}$, ${{\beta }_{i}}$  and introduce two one-parameter families of series $g_{\alpha }^{\left( t \right)}\left( x \right)$ and $g_{\beta }^{\left( t \right)}\left( x \right)$ associated with these expansions.

\end{abstract}

Matrices that we will consider correspond to operators in the space of formal power series over the field of real or complex numbers. Based on this, we associate the rows and columns of  matrices with the generating functions of their elements, i.e., formal power series. $n$th coefficient of the series $a\left( x \right)$ and $\left( n,m \right)$th element of the matrix$A$ will be denoted respectively $\left[ {{x}^{n}} \right]a\left( x \right)$, ${{\left( A \right)}_{n,m}}$. 

Matrix $\left( f\left( x \right),g\left( x \right) \right)$, $n$th column of which has the generating function $f\left( x \right){{g}^{n}}\left( x \right)$, ${{g}_{0}}=0$, is called Riordan matrix (Riordan array)  [1]. It is the product of two matrices that correspond to the operators of multiplication and composition of series:
$$\left( f\left( x \right),g\left( x \right) \right)=\left( f\left( x \right),x \right)\left( 1,g\left( x \right) \right),$$
$$\left( f\left( x \right),x \right)a\left( x \right)=f\left( x \right)a\left( x \right), \qquad\left( 1,g\left( x \right) \right)a\left( x \right)=a\left( g\left( x \right) \right),$$
$$\left( f\left( x \right),g\left( x \right) \right)\left( b\left( x \right),a\left( x \right) \right)=\left( f\left( x \right)b\left( g\left( x \right) \right),a\left( g\left( x \right) \right) \right).$$
Matrices $\left( f\left( x \right),g\left( x \right) \right)$, ${{f}_{0}}\ne 0$, ${{g}_{1}}\ne 0$, form a group called the Riordan group.
Matrices of the form $\left( f\left( x \right),x \right)$ form a subgroup called the Appell subgroup, matrices of the form $\left( 1,g\left( x \right) \right)$ form a subgroup called the associated subgroup, or the Lagrange subgroup.
 
 Matrices of the form
$$\left( a\left( x \right),x \right)+\left( xb\left( x \right),x \right)D,$$
where $D$ is the matrix of the differentiation operator: $D{{x}^{n}}=n{{x}^{n-1}}$, and on the coefficients of the series $a\left( x \right)$, $b\left( x \right)$ conditions are not imposed, form the Lie algebra of the Riordan group [2], [3]. 
Subalgebra of the matrices $\left( a\left( x \right),x \right)$ corresponds to the Appel subgroup, subalgebra of the matrices  $\left( xb\left( x \right),x \right)D$ corresponds to the Lagrange subgroup.
Lie algebra of the group of formal power series isomorphic to the Lagrange subgroup was considered in [4].

We restrict ourselves to the condition ${{g}_{1}}=1$. Then
$${{\left( \frac{g\left( x \right)}{x},x \right)}^{z}}=\sum\limits_{n=0}^{\infty }{\frac{{{z}^{n}}}{n!}}{{\left( \log \left( \frac{g\left( x \right)}{x},x \right) \right)}^{n}}, \quad{{\left( 1,g\left( x \right) \right)}^{t}}=\sum\limits_{n=0}^{\infty }{\frac{{{t}^{n}}}{n!}}{{\left( \log \left( 1,g\left( x \right) \right) \right)}^{n}},$$
where 
$$\log \left( \frac{g\left( x \right)}{x},x \right)=\left( \log \frac{g\left( x \right)}{x},x \right), \qquad\left( \log \left( 1,g\left( x \right) \right) \right)=\left( \omega \left( x \right),x \right)D,$$
$$\omega \left( x \right)=\sum\limits_{n=1}^{\infty }{{{\omega }_{n}}}{{x}^{n+1}}, \qquad\omega \left( g\left( x \right) \right)=\omega \left( x \right){g}'\left( x \right).$$
Denote ${{\left( 1,g\left( x \right) \right)}^{t}}=\left( 1,{{g}^{\left( t \right)}}\left( x \right) \right)$. Then
$${{\left( {g\left( x \right)}/{x}\; \right)}^{z}}=\sum\limits_{n=0}^{\infty }{\frac{{{z}^{n}}}{n!}}{{\left( \log \left( {g\left( x \right)}/{x}\; \right) \right)}^{n}}=\sum\limits_{n=0}^{\infty }{{{s}_{n}}}\left( z \right){{x}^{n}},$$
$${{g}^{\left( t \right)}}\left( x \right)=\sum\limits_{n=0}^{\infty }{\frac{{{t}^{n}}}{n!}}{{\omega }_{n}}\left( x \right)=\sum\limits_{n=1}^{\infty }{{{c}_{n}}}\left( t \right){{x}^{n}}, \quad{{\omega }_{0}}\left( x \right)=x, \quad{{\omega }_{n}}\left( x \right)=\omega \left( x \right){{{\omega }'}_{n-1}}\left( x \right).$$
Polynomials ${{s}_{n}}\left( z \right)$ are called convolution polynomials [5]. Explicit form of these polynomials is:
$${{s}_{n}}\left( z \right)=\sum\limits_{m=0}^{n}{{{z}^{m}}\sum\limits_{n,m}{\frac{\lambda _{1}^{{{m}_{1}}}\lambda _{2}^{{{m}_{2}}}...\text{ }\lambda _{n}^{{{m}_{n}}}}{{{m}_{1}}!{{m}_{2}}!\text{ }...\text{ }{{m}_{n}}!}}}\text{ },$$
where ${{\lambda }_{i}}=\left[ {{x}^{i}} \right]\left( \log \left( {g\left( x \right)}/{x}\; \right) \right)$ and  the summation of the coefficient of ${{z}^{m}}$ is over all partitions $\sum\nolimits_{i=1}^{n}{i{{m}_{i}}}=n$, $\sum\nolimits_{i=1}^{n}{{{m}_{i}}}=m$. Polynomials ${{c}_{n}}\left( x \right)$ will be called composition polynomials. Polynomials ${{s}_{n}}\left( t{{\omega }_{i}},z \right)$ such that
$${{\left( {{{g}^{\left( t \right)}}\left( x \right)}/{x}\; \right)}^{z}}=\sum\limits_{n=0}^{\infty }{{{s}_{n}}}\left( t{{\omega }_{i}},z \right){{x}^{n}}, \quad{{s}_{n}}\left( {{\omega }_{i}},z \right)={{s}_{n}}\left( z \right), \quad{{s}_{n}}\left( t{{\omega }_{i}},1 \right)={{c}_{n+1}}\left( t \right),$$
will be called the composition-convolution polynomials. Composition polynomials together with convolution polynomials were considered in [6]. It was shown that the composition-convolution polynomials have the form
$${{s}_{n}}\left( t{{\omega }_{i}},z \right)=\sum\limits_{m=0}^{n}{\frac{{{t}^{m}}}{m!}}\sum\limits_{n,m}{z\left( z+{{i}_{1}} \right)\left( z+{{i}_{1}}+{{i}_{2}} \right)...\left( z+{{i}_{1}}+{{i}_{2}}+...+{{i}_{m-1}} \right){{\omega }_{{{i}_{1}}}}{{\omega }_{{{i}_{2}}}}...{{\omega }_{{{i}_{m}}}}},$$
where the summation of the coefficient of ${{{t}^{m}}}/{m!}\;$ is over all compositions $n={{i}_{1}}+{{i}_{2}}+...+{{i}_{m}}$. Earlier this formula was obtained in [7].

In this note, we give two formulas similar to the formula for polynomials ${{s}_{n}}\left( t{{\omega }_{i}},z \right)$ and probably related to it. The starting point for further constructions is the formula
$$\left( \alpha {{x}^{k+1}},x \right)D=\log \left( 1,\frac{x}{{{\left( 1-k\alpha {{x}^{k}} \right)}^{{1}/{k}\;}}} \right), \qquad k>0;$$
$$\frac{1}{{{\left( 1-k\alpha {{x}^{k}} \right)}^{{m}/{k}\;}}}=\sum\limits_{n=0}^{\infty }{\frac{m\left( m+k \right)\left( m+2k \right)...\left( m+\left( n-1 \right)k \right)}{n!}{{\alpha }^{n}}{{x}^{nk}}}.$$
We denote the matrix whose $m$th column has the generating function ${{x}^{m}}{{\left( 1-k\alpha {{x}^{k}} \right)}^{{-m}/{k}\;}}$ by $P_{k}^{\alpha }$. Matrix $\left( 1,g\left( x \right) \right)$ can be decomposed into an infinite product of the matrices $P_{k}^{\alpha }$ with suitable exponents in two ways  –  to left-handed and right-handed products with respect to the matrix $P_{1}^{{{\alpha }_{1}}={{\beta }_{1}}}$:
$$\left( 1,g\left( x \right) \right)=...P_{k}^{{{\alpha }_{k}}}...P_{2}^{{{\alpha }_{2}}}P_{1}^{{{\alpha }_{1}}}=P_{1}^{{{\beta }_{1}}}P_{2}^{{{\beta }_{2}}}...P_{k}^{{{\beta }_{k}}}....$$
The $\left( n,m \right)$th element of the matrix $\left( 1,g\left( x \right) \right)$, expressed through the elements of the matrices $P_{1}^{{{\alpha }_{1}}},...,P_{n-1}^{{{\alpha }_{n-1}}}$ or $P_{1}^{{{\beta }_{1}}},...,P_{n-1}^{{{\beta }_{n-1}}}$, can be represented as a polynomial in $m$ of degree $n-m$. Replacing  $m$ with $z$, we obtain two formulas for polynomials ${{s}_{n}}\left( z \right)=\left[ {{x}^{n}} \right]{{\left( {g\left( x \right)}/{x}\; \right)}^{z}}$ expressed in terms of expansion coefficients ${{\alpha }_{i}}$, ${{\beta }_{i}}$:
$${{s}_{1}}\left( {{\alpha }_{i}},z \right)=z{{\alpha }_{1}}, \qquad{{s}_{2}}\left( {{\alpha }_{i}},z \right)=z{{\alpha }_{2}}+\left( \begin{matrix}
   z+1  \\
   2  \\
\end{matrix} \right)\alpha _{1}^{2},$$
$${{s}_{3}}\left( {{\alpha }_{i}},z \right)=z{{\alpha }_{3}}+z\left( z+1 \right){{\alpha }_{1}}{{\alpha }_{2}}+\left( \begin{matrix}
   z+2  \\
   3  \\
\end{matrix} \right)\alpha _{1}^{3},$$
$${{s}_{4}}\left( {{\alpha }_{i}},z \right)=z{{\alpha }_{4}}+z\left( z+1 \right){{\alpha }_{1}}{{\alpha }_{3}}+z\left( z+2 \right)\frac{\alpha _{2}^{2}}{2}+
z\left( z+1 \right)\left( z+2 \right)\frac{\alpha _{1}^{2}{{\alpha }_{2}}}{2}+\left( \begin{matrix}
   z+3  \\
   4  \\
\end{matrix} \right)\alpha _{1}^{4};$$

$${{s}_{1}}\left( {{\beta }_{i}},z \right)=z{{\beta }_{1}},  \qquad{{s}_{2}}\left( {{\beta }_{i}},z \right)=z{{\beta }_{2}}+\left( \begin{matrix}
   z+1  \\
   2  \\
\end{matrix} \right)\beta _{1}^{2},$$
$${{s}_{3}}\left( {{\beta }_{i}},z \right)=z{{\beta }_{3}}+z\left( z+2 \right){{\beta }_{1}}{{\beta }_{2}}+\left( \begin{matrix}
   z+2  \\
   3  \\
\end{matrix} \right)\beta _{1}^{3},$$
$${{s}_{4}}\left( {{\beta }_{i}},z \right)=z{{\beta }_{4}}+z\left( z+3 \right){{\beta }_{1}}{{\beta }_{3}}+z\left( z+2 \right)\frac{\beta _{2}^{2}}{2}+
z\left( z+2 \right)\left( z+3 \right)\frac{\beta _{1}^{2}{{\beta }_{2}}}{2}+\left( \begin{matrix}
   z+3  \\
   4  \\
\end{matrix} \right)\beta _{1}^{4}.$$
{\bfseries Theorem 1.} \emph{. Polynomials ${{s}_{n}}\left( {{\alpha }_{i}},z \right)$, ${{s}_{n}}\left( {{\beta }_{i}},z \right)$  have the form
$${{s}_{n}}\left( {{\alpha }_{i}},z \right)=\sum\limits_{n}{z\left( z+{{i}_{1}} \right)\left( z+{{i}_{1}}+{{i}_{2}} \right)...\left( z+{{i}_{1}}+{{i}_{2}}+...+{{i}_{m-1}} \right)\frac{\alpha _{1}^{{{m}_{1}}}\alpha _{2}^{{{m}_{2}}}...\text{ }\alpha _{n}^{{{m}_{n}}}}{{{m}_{1}}!{{m}_{2}}!...{{m}_{n}}!}}\text{ },$$
$${{s}_{n}}\left( {{\beta }_{i}},z \right)=\sum\limits_{n}{z\left( z+{{i}_{1}} \right)\left( z+{{i}_{1}}+{{i}_{2}} \right)...\left( z+{{i}_{1}}+{{i}_{2}}+...+{{i}_{m-1}} \right)\frac{\beta _{1}^{{{m}_{1}}}\beta _{2}^{{{m}_{2}}}...\text{ }\beta _{n}^{{{m}_{n}}}}{{{m}_{1}}!{{m}_{2}}!...{{m}_{n}}!}}\text{ },$$
where the summation is over all partitions $n=\sum\nolimits_{i=1}^{n}{i{{m}_{i}}}={{i}_{1}}+{{i}_{2}}+...+{{i}_{m}}$, $m=\sum\nolimits_{i=1}^{n}{{{m}_{i}}}$, but for the first formula ${{i}_{k}}\le {{i}_{k+1}}$, for the second ${{i}_{k}}\ge {{i}_{k+1}}$.
}\\
{\bfseries Proof.} Let ${{A}_{k}}$ is the matrix whose $i$th column has the generating function ${{x}^{i}}\sum\nolimits_{n=0}^{\infty }{_{\left( k \right)}a_{n}^{\left( m+i \right)}}{{x}^{kn}}$, $m=1$, $2$, …. Then
$${{\left( {{A}_{2}}{{A}_{1}} \right)}_{n,0}}=\sum\limits_{{{m}_{1}}+2{{m}_{2}}=n}{_{\left( 1 \right)}a_{{{m}_{1}}}^{\left( m \right)}}{}_{\left( 2 \right)}a_{{{m}_{2}}}^{\left( m+{{m}_{1}} \right)},$$
$${{\left( {{A}_{3}}{{A}_{2}}{{A}_{1}} \right)}_{n,0}}=\sum\limits_{{{m}_{1}}+2{{m}_{2}}+3{{m}_{3}}=n}{_{\left( 1 \right)}a_{{{m}_{1}}}^{\left( m \right)}}{}_{\left( 2 \right)}a_{{{m}_{2}}}^{\left( m+{{m}_{1}} \right)}{}_{\left( 3 \right)}a_{{{m}_{3}}}^{\left( m+{{m}_{1}}+2{{m}_{2}} \right)},$$
$${{\left( ...{{A}_{k}}...{{A}_{2}}{{A}_{1}} \right)}_{n,0}}=\sum\limits_{{{m}_{1}}+2{{m}_{2}}+...+n{{m}_{n}}=n}{_{\left( 1 \right)}a_{{{m}_{1}}}^{\left( m \right)}}{}_{\left( 2 \right)}a_{{{m}_{2}}}^{\left( m+{{m}_{1}} \right)}...{}_{\left( n \right)}a_{{{m}_{n}}}^{\left( m+{{m}_{1}}+...+\left( n-1 \right){{m}_{n-1}} \right)};$$
$${{\left( {{A}_{1}}{{A}_{2}} \right)}_{n,0}}=\sum\limits_{{{m}_{1}}+2{{m}_{2}}=n}{_{\left( 2 \right)}a_{{{m}_{2}}}^{\left( m \right)}}{}_{\left( 1 \right)}a_{{{m}_{1}}}^{\left( m+2{{m}_{2}} \right)},$$
$${{\left( {{A}_{1}}{{A}_{2}}{{A}_{3}} \right)}_{n,0}}=\sum\limits_{{{m}_{1}}+2{{m}_{2}}+3{{m}_{3}}=n}{_{\left( 3 \right)}a_{{{m}_{3}}}^{\left( m \right)}}{}_{\left( 2 \right)}a_{{{m}_{2}}}^{\left( m+3{{m}_{3}} \right)}{}_{\left( 1 \right)}a_{{{m}_{1}}}^{\left( m+3{{m}_{3}}+2{{m}_{2}} \right)},$$
$${{\left( {{A}_{1}}{{A}_{2}}...{{A}_{k}}... \right)}_{n,0}}=\sum\limits_{{{m}_{1}}+2{{m}_{2}}+...+n{{m}_{n}}=n}{_{\left( n \right)}a_{{{m}_{n}}}^{\left( m \right)}}{}_{\left( n-1 \right)}a_{{{m}_{n-1}}}^{\left( m+n{{m}_{n}} \right)}...{}_{\left( 1 \right)}a_{{{m}_{1}}}^{\left( m+...+2{{m}_{2}} \right)}.$$
If 
$$_{\left( k \right)}a_{0}^{\left( m \right)}=1,  \qquad_{\left( k \right)}a_{{{m}_{k}}}^{\left( m \right)}=m\left( m+k \right)\left( m+2k \right)...\left( m+\left( {{m}_{k}}-1 \right)k \right)\frac{\alpha _{k}^{{{m}_{k}}}}{{{m}_{k}}!},$$
then
$${{\left( ...{{A}_{k}}...{{A}_{2}}{{A}_{1}} \right)}_{n,0}}=\sum\limits_{{{m}_{1}}+2{{m}_{2}}+...+n{{m}_{n}}=n}{m\left( m+1 \right)}\left( m+2 \right)...\left( m+{{m}_{1}}-1 \right)\times $$
$$\times \left( m+{{m}_{1}} \right)\left( m+{{m}_{1}}+2 \right)...\left( m+{{m}_{1}}+\left( {{m}_{2}}-1 \right)2 \right)\times ...$$
$$...\times \left( m+{{m}_{1}}+2{{m}_{2}}+...+\left( n-1 \right){{m}_{n-1}}+\left( {{m}_{n}}-1 \right)n \right)\frac{\alpha _{1}^{{{m}_{1}}}\alpha _{2}^{{{m}_{2}}}...\text{ }\alpha _{n}^{{{m}_{n}}}}{{{m}_{1}}!{{m}_{2}}!...{{m}_{n}}!}={{s}_{n}}\left( {{\alpha }_{i}},m \right),$$
and hence ${{\left( ...P_{k}^{{{\alpha }_{k}}}...P_{2}^{{{\alpha }_{2}}}P_{1}^{{{\alpha }_{1}}} \right)}_{n-m,m}}={{s}_{n-m}}\left( {{\alpha }_{i}},m \right)$. Similarly, if
$$_{\left( k \right)}a_{0}^{\left( m \right)}=1,  \qquad_{\left( k \right)}a_{{{m}_{k}}}^{\left( m \right)}=m\left( m+k \right)\left( m+2k \right)...\left( m+\left( {{m}_{k}}-1 \right)k \right)\frac{\beta _{k}^{{{m}_{k}}}}{{{m}_{k}}!},$$
then
$${{\left( {{A}_{1}}{{A}_{2}}...{{A}_{k}}... \right)}_{n,0}}=\sum\limits_{{{m}_{1}}+2{{m}_{2}}+...+n{{m}_{n}}=n}{m\left( m+n \right)}...\left( m+\left( {{m}_{n}}-1 \right)n \right)\times $$
$$\times \left( m+n{{m}_{n}} \right)\left( m+n{{m}_{n}}+\left( n-1 \right) \right)...\left( m+n{{m}_{n}}+\left( {{m}_{n-1}}-1 \right)\left( n-1 \right) \right)\times ...$$
$$...\times \left( m+n{{m}_{n}}+...+2{{m}_{2}}+\left( {{m}_{1}}-1 \right) \right)\frac{\beta _{1}^{{{m}_{1}}}\beta _{2}^{{{m}_{2}}}...\text{ }\beta _{n}^{{{m}_{n}}}}{{{m}_{1}}!{{m}_{2}}!...{{m}_{n}}!}={{s}_{n}}\left( {{\beta }_{i}},m \right),$$   
$${{\left( P_{1}^{{{\beta }_{1}}}P_{2}^{{{\beta }_{2}}}...P_{k}^{{{\beta }_{k}}}... \right)}_{n-m,m}}={{s}_{n-m}}\left( {{\beta }_{i}},m \right).   \qquad      \square $$

Note that since
 $${{\left( 1,g\left( x \right) \right)}^{-1}}=...P_{k}^{-{{\beta }_{k}}}...P_{2}^{-{{\beta }_{2}}}P_{1}^{-{{\beta }_{1}}}=P_{1}^{-{{\alpha }_{1}}}P_{2}^{-{{\alpha }_{2}}}...P_{k}^{-{{\alpha }_{k}}}...,$$
 then in accordance with the Lagrange inversion theorem
$$\frac{z}{z+n}{{s}_{n}}\left( {{\alpha }_{i}},-z-n \right)={{s}_{n}}\left( {{\beta }_{i}},z \right), \qquad{{\beta }_{i}}=-{{\alpha }_{i}};$$ 
$$\frac{z}{z+n}{{s}_{n}}\left( {{\beta }_{i}},-z-n \right)={{s}_{n}}\left( {{\alpha }_{i}},z \right), \qquad{{\alpha }_{i}}=-{{\beta }_{i}}.$$  

Denote $\alpha \left( x \right)=\sum\nolimits_{n=1}^{\infty }{{{\alpha }_{n}}{{x}^{n+1}}}$, $\beta \left( x \right)=\sum\nolimits_{n=1}^{\infty }{{{\beta }_{n}}{{x}^{n+1}}}$. Note that if the matrix   $\left( 1,g\left( x \right) \right)$ is a pseudo-involution, i.e., ${{g}^{\left( -1 \right)}}\left( x \right)=-g\left( -x \right)$, then $\omega \left( -x \right)=\omega \left( x \right)$, $\alpha \left( -x \right)=\beta \left( x \right)$.

We introduce two one-parameter families of series $g_{\alpha }^{\left( t \right)}\left( x \right)$ and $g_{\beta }^{\left( t \right)}\left( x \right)$:
$${{\left( {g_{\alpha }^{\left( t \right)}\left( x \right)}/{x}\; \right)}^{z}}=\sum\limits_{n=0}^{\infty }{{{s}_{n}}}\left( t{{\alpha }_{i}},z \right),  \qquad{{\left( {g_{\beta }^{\left( t \right)}\left( x \right)}/{x}\; \right)}^{z}}=\sum\limits_{n=0}^{\infty }{{{s}_{n}}}\left( t{{\beta }_{i}},z \right);$$
$${{s}_{1}}\left( t{{\alpha }_{i}},1 \right)={{\alpha }_{1}}t, \qquad{{s}_{2}}\left( t{{\alpha }_{i}},1 \right)={{\alpha }_{2}}t+\alpha _{1}^{2}{{t}^{2}}, \qquad{{s}_{3}}\left( t{{\alpha }_{i}},1 \right)={{\alpha }_{3}}t+2{{\alpha }_{1}}{{\alpha }_{2}}{{t}^{2}}+\alpha _{1}^{3}{{t}^{3}},$$
$${{s}_{4}}\left( t{{\alpha }_{i}},1 \right)={{\alpha }_{4}}t+\left( 2{{\alpha }_{1}}{{\alpha }_{3}}+\frac{3}{2}\alpha _{2}^{2} \right){{t}^{2}}+3\alpha _{1}^{2}{{\alpha }_{2}}{{t}^{3}}+\alpha _{1}^{4}{{t}^{4}},$$
$${{s}_{5}}\left( t{{\alpha }_{i}},1 \right)={{\alpha }_{5}}t+\left( 2{{\alpha }_{1}}{{\alpha }_{4}}+3{{\alpha }_{2}}{{\alpha }_{3}} \right){{t}^{2}}+\left( 3\alpha _{1}^{2}{{\alpha }_{3}}+3{{\alpha }_{1}}\alpha _{2}^{2} \right){{t}^{3}}+4\alpha _{1}^{3}{{\alpha }_{2}}{{t}^{4}}+\alpha _{1}^{5}{{t}^{5}};$$
$${{s}_{1}}\left( t{{\beta }_{i}},1 \right)={{\beta }_{1}}t, \qquad{{s}_{2}}\left( t{{\beta }_{i}},1 \right)={{\beta }_{2}}t+\beta _{1}^{2}{{t}^{2}}, \qquad{{s}_{3}}\left( t{{\beta }_{i}},1 \right)={{\beta }_{3}}t+3{{\beta }_{1}}{{\beta }_{2}}{{t}^{2}}+\beta _{1}^{3}{{t}^{3}},$$
$${{s}_{4}}\left( t{{\beta }_{i}},1 \right)={{\beta }_{4}}t+\left( 4{{\beta }_{1}}{{\beta }_{3}}+\frac{3}{2}\alpha _{2}^{2} \right){{t}^{2}}+6\beta _{1}^{2}{{\beta }_{2}}{{t}^{3}}+\beta _{1}^{4}{{t}^{4}},$$
$${{s}_{5}}\left( t{{\beta }_{i}},1 \right)={{\beta }_{5}}t+\left( 5{{\beta }_{1}}{{\beta }_{4}}+4{{\beta }_{2}}{{\beta }_{3}} \right){{t}^{2}}+\left( 10\beta _{1}^{2}{{\beta }_{3}}+\frac{15}{2}{{\beta }_{1}}\beta _{2}^{2} \right){{t}^{3}}+10\beta _{1}^{3}{{\beta }_{2}}{{t}^{4}}+\beta _{1}^{5}{{t}^{5}}.$$
Denote ${{\left( 1,g\left( x \right) \right)}^{-1}}=\left( 1,\bar{g}\left( x \right) \right)$. Then ${{\bar{\alpha }}_{i}}=-{{\beta }_{i}}$, ${{\bar{\beta }}_{i}}=-{{\alpha }_{i}}$, 
$${{\left( 1,g_{\alpha }^{\left( t \right)}\left( x \right) \right)}^{-1}}=\left( 1,\bar{g}_{\beta }^{\left( t \right)}\left( x \right) \right),  \qquad{{\left( 1,g_{\beta }^{\left( t \right)}\left( x \right) \right)}^{-1}}=\left( 1,\bar{g}_{\alpha }^{\left( t \right)}\left( x \right) \right).$$

The following two propositions are given without proofs.\\
{\bfseries Proposition 1.} The following formula is true
$$\left( 1,g_{\beta }^{\left( -1 \right)}\left( x \right) \right)\left( 1,g_{\alpha }^{\left( -1 \right)}\left( x \right) \right)=\left( 1,{{g}^{\left( -2 \right)}}\left( x \right) \right),$$
or
$$\left( P_{1}^{-{{\beta }_{1}}}P_{2}^{-{{\beta }_{2}}}...P_{k}^{-{{\beta }_{k}}}... \right)\left( ...P_{k}^{-{{\alpha }_{k}}}...P_{2}^{-{{\alpha }_{2}}}P_{1}^{-{{\alpha }_{1}}} \right)=$$
$$=\left( P_{1}^{-{{\alpha }_{1}}}P_{2}^{-{{\alpha }_{2}}}... \right)\left( ...P_{2}^{-{{\beta }_{2}}}P_{1}^{-{{\beta }_{1}}} \right)=\left( ...P_{2}^{-{{\beta }_{2}}}P_{1}^{-{{\beta }_{1}}} \right)\left( P_{1}^{-{{\alpha }_{1}}}P_{2}^{-{{\alpha }_{2}}}... \right).$$
Respectively,
$$\left( 1,\bar{g}_{\beta }^{\left( -1 \right)}\left( x \right) \right)\left( 1,\bar{g}_{\alpha }^{\left( -1 \right)}\left( x \right) \right)=\left( 1,{{g}^{\left( 2 \right)}}\left( x \right) \right),$$
or
$$\left( P_{1}^{{{\alpha }_{1}}}P_{2}^{{{\alpha }_{2}}}... \right)\left( ...P_{2}^{{{\beta }_{2}}}P_{1}^{{{\beta }_{1}}} \right)=\left( P_{1}^{{{\beta }_{1}}}P_{2}^{{{\beta }_{2}}}... \right)\left( ...P_{2}^{{{\alpha }_{2}}}P_{1}^{{{\alpha }_{1}}} \right)=\left( ...P_{2}^{{{\alpha }_{2}}}P_{1}^{{{\alpha }_{1}}} \right)\left( P_{1}^{{{\beta }_{1}}}P_{2}^{{{\beta }_{2}}}... \right).$$\\
{\bfseries Proposition 2.} The following formula is true.
$$\left( 1,g_{\alpha }^{\left( t \right)}\left( x \right) \right)\left( 1,g_{\beta }^{\left( 1-t \right)}\left( x \right) \right)=\left( 1,g\left( x \right) \right),$$
or
$$\left( 1,g\left( x \right) \right)=...P_{k}^{t{{\alpha }_{k}}}...P_{2}^{t{{\alpha }_{2}}}P_{1}^{{{\alpha }_{1}}={{\beta }_{1}}}P_{2}^{\left( 1-t \right){{\beta }_{2}}}...P_{k}^{\left( 1-t \right){{\beta }_{k}}}....$$

Let ${{D}_{t}}$ is the differentiation operator with respect to $t$. Then
$${{D}_{t}}g_{\alpha }^{\left( t \right)}\left( x \right){{|}_{t=0}}=\alpha \left( x \right), \qquad{{D}_{t}}g_{\beta }^{\left( t \right)}\left( x \right){{|}_{t=0}}=\beta \left( x \right),$$
$${{D}_{t}}\bar{g}_{\alpha }^{\left( t \right)}\left( x \right){{|}_{t=0}}=-\beta \left( x \right), \qquad{{D}_{t}}\bar{g}_{\beta }^{\left( t \right)}\left( x \right){{|}_{t=0}}=-\alpha \left( x \right).$$
Denote 
$${{D}_{t}}g_{\alpha }^{\left( t \right)}\left( x \right){{|}_{t=1}}={{f}_{\alpha }}\left( x \right), \qquad{{D}_{t}}g_{\beta }^{\left( t \right)}\left( x \right){{|}_{t=1}}={{f}_{\beta }}\left( x \right).$$
Then
$${{D}_{t}}g_{\alpha }^{\left( 1-t \right)}\left( x \right){{|}_{t=0}}={{D}_{t}}\bar{g}_{\alpha }^{\left( t \right)}\left( g\left( x \right) \right){{|}_{t=0}},  \qquad{{D}_{t}}g_{\beta }^{\left( 1-t \right)}\left( x \right){{|}_{t=0}}={{D}_{t}}g\left( \bar{g}_{\beta }^{\left( t \right)}\left( x \right) \right){{|}_{t=0}},$$
or, since${{D}_{t}}{{\left( \bar{g}_{\beta }^{\left( t \right)}\left( x \right) \right)}^{n}}{{|}_{t=0}}=-n{{x}^{n-1}}\alpha \left( x \right)$,
$${{f}_{\alpha }}\left( x \right)=\beta \left( g\left( x \right) \right),  \qquad{{f}_{\beta }}\left( x \right)=\alpha \left( x \right){g}'\left( x \right).$$

E-mail: {evgeniy\symbol{"5F}burlachenko@list.ru}
\end{document}